\documentclass{amsart}
\usepackage{amsfonts}

\begin{document}

\title{ON THE EXISTENCE OF CERTAIN METRICS WITH NON PARALLEL RICCI TENSOR }

\begin{center}
\author{A. Raouf  Chouikha  }

\thanks{Universite Paris 13 LAGA UMR 7539, Villetaneuse 93430}
 
\begin{abstract}
A. Derdzinki [D] gave examples of Riemannian metrics with harmonic curvature and non parallel Ricci tensor on some compact  manifolds \ $(M,g]$\ . We examine their existence as well as their number wich naturally depends on the geometry of the manifolds.
\end{abstract}
\maketitle

\end{center}

\section{Introduction}

 Let \ $(M,g)$\ be a Riemannian manifold of dimension $n \geq 3$.\ $M$\  is said to have harmonic curvature if the divergence of its curvature tensor \ ${\mathcal R}$\ vanishes (in local coordonates : $\nabla^i {\mathcal R}_{ijkl} = 0$). That means the Ricci tensor \ $r$\ is a codazzi tensor \ ($\nabla^i {\mathcal R}_{ijkl} = \nabla_k r_{hj} - \nabla_j r_{hk} = 0$).\ On other words, in the compact case of the manifold the Riemannian connection is a Yang-Mills potential in the tangent bundle.\\
Answering to the question on the parallelism of the Ricci tensor of the Riemannian metrics, A. Derdzinski gave examples of compact manifolds with with harmonic curvature but non parallel Ricci tensor: \ $\delta {\mathcal R} = 0$\ and \ $\nabla r \neq 0$.\ Moreover, he obtains some classification results, [D]. \\
The corresponding manifolds are bundles with fibres \ $N$
\ over the circle \
$S^1$ \ (parametrized by arc length \ t\ and of lenght \ $T = \int_{S^1} dt$) \ equipped with the warped metrics
\quad $dt^2 + h^{4/n}(t)g_0$ \quad on the product  \quad
$ S^1\times N $ .\quad Here, \ $(N,g_0)$ \ is an Einstein manifold of dimension \ $n \geq 3$, \  and the
function \ $h(t)$ \ on the prime factor is a periodic solution of the ODE, established by 
Derdzinski 
\begin{equation}
 h'' - {nR\over{4(n-1)}}h^{1-4/n} = -{n\over 4}C h \qquad \mbox{for some constant}
\qquad C > 0.
\end{equation}
 
 This function must be
non constant, otherwise the corresponding metric has a parallel Ricci tensor. \\
The goal of this paper is to study the existence and  the number of such metric which naturally must depend on the geometry of \ $S^1$ \ and \ $N$. In particular, we prove there is a positive bound \ $T_0$\ depending on the \ $S^1$ \ lenght  such that when $$T < T_0$$ the above equation may have only constant solutions, i.e.
$$h(t) \equiv \alpha = \bigg({R\over{4(n-1) C}}\bigg)^{4/n}.$$ The corresponding Ricci tensor cannot be non parallel.\\ 

\section{Metrics with harmonic curvature}
Let \ $(M,g)$\ be a Riemannian manifold with \ $n = dim M \geq 3. \ {\mathcal R}$\ is its curvature tensor, \ $r$\ its Ricci tensor \ $W$ \ its Weyl conformal tensor and \ $R$ \ its scalar curvature. According to the second Bianchi identity \ $d {\mathcal R} = 0$ \ (in local coordonates \ $\nabla_q {\mathcal R}_{ijkl} + \nabla_i {\mathcal R}_{jqkl} + \nabla_j {\mathcal R}_{qikl} = 0,$)\ we get the following relations
$$ \delta {\mathcal R} = - d r \qquad i.e. \quad \nabla^i {\mathcal R}_{ijkl} = \nabla_k r_{hj} - \nabla_j r_{hk}$$
$$( n - 2) \delta W = - (n - 3) d[r - \frac {R g}{2n - 2}] \qquad {\it and}\quad 2 \delta r = - d R.$$
The sign conventions are such that \ $ r_{ij} = {\mathcal R}_{ilj}^l, \quad R = g^{ij} r_{ij}$.\\
$d$\ is the exterior differentiation and \ $\delta$ is its formal adjoint, viewed them as differential forms on \ $M$.\\
$(M,g)$\ has harmonic curvature if \ $\delta {\mathcal R} = 0.$\\
A Codazzi tensor \ $C$\ on \ $(M,g)$\ is a symetric \ $(0,2)$\ tensor field on \ $M$\ verifying the Codazzi equation \ $d C = 0\quad i.e. \ \nabla_j C_{ik} = \nabla_k C_{ij}.$\\
Classical properties of metrics with harmonic curvature are summarized in the following lemma, [B], [S]

\bigskip
{\bf Lemma 1}\quad {\it Let \ $(M,g]$\ be a Riemannian manifold of dimension $n \geq 3$. The following holds \\
(1) \ if \ $n = 3,\  (M,g)$ has harmonic curvature \ $\delta {\mathcal R} = 0$ \ if and only if it is conformally flat ($W \equiv 0$) and has constant scalar curvature \ $R = Cte$.\\
(2) \ If \ $n \geq 4,\ \delta W = 0 $\ and \ $M$ \ has constant scalar curvature, then \ $\delta {\mathcal R} = 0$ .\\
(3) \ $(M,g)$\ has harmonic curvature if and only if its Ricci tensor is codazzi \ ($d r = 0$.\\
(4) \ If \ $(M,g)$\ is a Riemannian product, then it has harmonic curvature if and only if any factor manifolds has harmonic curvature.}

\bigskip

 More generally, Derdzinski established a classification of the compact n-dimensional Riemannian
manifolds\quad $(M_n,g), \ n \geq 3$, \quad with harmonic curvature. If the Ricci tensor \ $Ric(g)$ \
is not parallel and has less than three distinct eigenvalues at each point, then \ $(M,g)$ \ is
covered isometrically by a manifold  $$(S^1(T) \times N,dt^2 + h^{4/n}(t)g_0),$$ where
the non constant positive periodic solutions \ $h$ \ verify the equation (1). Here \ $(N,g_0)$\ is a (n-1)-
dimensional Einstein manifold with positive (constant) scalar curvature.\\

 \section{
 On the existence of Derdzinski metrics}
\subsection{Analysis of the Derdzinski equation}
The function $h$ defined by the above warped product \ $(S^1(T) \times N,dt^2 + h^{4/n}(t)g_0)$\ is a solution of 
$$(1) \qquad h'' - {nR\over{4(n-1)}}h^{1-4/n} = -{n\over 4}C h \qquad \mbox{for some constant}
\qquad C > 0$$
so that its curvature is harmonic and its Ricci tensor is non parallel \ $ Dr
\neq 0$.\\
In fact, this may be deduced from the following result (Lemma 1 in [D].

\bigskip

{\bf Lemma 2}\quad {\it Let $I$ be an interval of $\!R$, \ and  $q$ a $C^{\infty}$ function on $I$ such that $ e^q = h^{4/n}$.\ Let  $(N,g_0)$ be an $(n-1)$ dimensional Riemannian manifold, $r_0$ its Ricci tensor, and $R$ its scalar curvature with $n \geq 3$. Consider the warped product $I \times _{e^q} N,dt^2 + h^{4/n}(t)g_0)$\  and $r$ its Ricci tensor.\\
(1) \quad Given a local product chart $t = x^0, x^1, x^2, .....,x^{n-1} $ for $I\times N$ with $g_{00} = 1, g_{0i} = 0$ and $g_{ij} = e^q g_{0_{ij}}.$ \ The components of the covariant derivative of $r$ are 
$$\nabla_0 r_{00} = - \frac {n-1}{2} [q''' + q' q''], \quad \nabla_0 r_{i0} =  \nabla_i r_{00} = 0,$$
$$\nabla_0 r_{ij} = - q' r_{0_{ij}} - \frac{1}{2} e^q [q''' + (n-1) q'q''] g_{0_{ij}}, \quad \nabla_i r_{0j} = - \frac{1}{2} q' r_{0_{ij}} - \frac{n-2}{4} e^q  q'q'' g_{0_{ij}}.$$
(2)\quad If $q$ is non constant, then the product $I \times _{e^q} N,dt^2 + h^{4/n}(t)g_0)$\ has harmonic curvature if and only if $(N,g_0)$ is an Einstein manifold and the positive function $h$ satisfies the ODE (1) on $I$ .}

\bigskip
Notice that the warped product manifold \ $(\!R \times M', dt^2 + h^{4/n}(t)g')$\ is complete analytic.
It thus appears from Lemma 2 that condition \ $\nabla r = 0$\ is equivalent to \  $q''' + q'q'' = 0$ \ and \ $(n-1)(n-2) q'' e^q + 2 R = 0$.\_
Moreover, [D] proved (his Theorem 1) that equation (1) posseses at least one non constant positive periodic solution if and only if $R > 0$ and $- {n\over 4}p = C > 0$.\\  
But, some care needs to be taken on the subject of the constant and the existence of a non trivial solution. This constant 
\ $C$\ must to verify some additional condition so that the conclusion holds.\\
More precisely, as we shall see below it does not exist any non constant periodic solution of (1)  with period $T$ if the constant \ $C$\ verifies
$$0 < C \leq \frac {4 \pi ^2}{T^2}.$$

\subsection{Existence conditions of the Derdzinski metrics}\ 
\newline More precisely, in considering the ODE point of view, we are able to analyse all the
solutions of equation (1).\\
So, this equation may have a non constant, positive periodic solution.\\
We shall prove the following

\bigskip

{\bf Theorem 1 }\qquad {\it Consider the warped product \quad $(S^1(T) \times N,dt^2 + h^{4/n}(t)g_0).$\quad 
The non constant positive periodic solution \ $h(t)$\ satisfies Equation (1) and \ $(N,g_0)$\
is a (n-1)-dimensional Einstein manifold with positive (constant) scalar curvature \ $R$\ . Then,
if the circle length \ $T$\ satisfies the following inequalities $$
2\pi{{(k-1)}\over{\sqrt C}} <T \leq {k\over{\sqrt C}}\ , \qquad \mbox{where\  k \ is an integer} > 1,$$ there exist at least \ $k$\ rotationnally invariant warped metrics \\ $dt^2 + h^{4/n}(t)g_0$ \ on the product manifold \ $S^1(T)
\times N.$
\\ Moreover, these metrics have an harmonic curvature. Their Ricci tensor are non parallel only if \ $
 T >{{2\pi}\over{\sqrt C}}$\ .\\
 Conversely, if the warped metric of the type \ $dt^2 + h^{4/n}(t)g_0$ \ on the manifold \ $S^1(T)
\times N$
 has harmonic curvature, then the function \ $h(t)$ \ on the circle satisfies an ODE (1).}\\

\bigskip

Condition \ $
 T > {{2\pi}\over{\sqrt C}}$\ implies the existence of a non constant solution \ $h(t)$\ .
Otherwise,  a trivial product have a parallel Ricci tensor.

\subsection{Proof of Theorem 1}\qquad 
 All periodic orbits \ $\gamma_c(t)$ \ of Equation (1) are surrounded by the homoclinic orbit \
$\gamma_{c_0}$ .\ The last one may be parametrized by \quad $(u_0(t),v_0(t)).$\\
Denote the coordinates of \ $\gamma_c(t)$ \ by \ $(u_c(t),v_c(t).$ \
 When the value \ $c$ \ satisfies the condition \ $1 < c < c_0$ \ , the
correspondant orbit is periodic ($c_0$\ corresponds to a periodic solution of null energy) .\\
The center is \ $(\alpha , 0$, \ where \ $\alpha  = \bigg({\displaystyle {R\over{4(n-1) C}}}\bigg)^{4/n}$.\\
 One may easily remark that two positive T-periodic solutions of
(1) having the same energy are translated, and thus give rise to equivalent
metrics on \
$(S^1(T)\times S^{n-1},g_0)$. \\  Note that the metric corresponding to the conformal factor
\quad
$u_0\ : \quad  \ g = {u_0}^{4\over{n-2}}g_0$ \quad is non complete. Then, it is not a pseudo-cylindric
metric; for this reason the constant \ $c$ \ cannot attain the critical value \ $c_0$ \ .\\
Equation (1) may be written under the following form
\begin{equation}
x'' + \phi (x) = 0
\end{equation}
where
$$\phi (x) = {n\over 4}C (x - \alpha ) - {nR\over{4(n-1)}}(x - \alpha )^{1-4/n} .$$
The period of the periodic solutions depends on the energy \ $T\equiv T(c)$\  with
$c$ the energy constant. It can be expressed 
$$T(c) = {\sqrt 2} \int_a^b \frac {du }{\sqrt{c - G(u)}}$$
where $G(u)$ is an integral of $\phi (u),$ with a nondegenerate relative minimum
at the origin. It verifies in addition  \ $ G(a) = G(b) = c$\ and \ $ a \leq \alpha \leq b$.  \\
So, \ $\phi (\alpha ) = 0 $\ and \ $\phi '(\alpha ) = {n\over 4}C > 0.$\ Hence, the origin is a center of Equation (2). That means in the neighbourhood of the trivial solution \ $h(t) \equiv \alpha$\ Equation (1) admits  a periodic solution.\\ 
We need the following result\\

 {\bf Lemma 3 }\qquad {\it Under the above hypothesis the  family of solutions \ $(T,u_T(t))$ \ of the ODE (2) 
(where \ $T$ \ is the minimal period) has bifurcation points on the values
\quad $(T_k,u_T{_k}(t))$ \ where \ $ T_k={{2\pi k}\over{\sqrt{n-2}}}$\ and \ $u_{T_k}\equiv\alpha$\ is a constant.  
In this family, there is a curve of non trivial solutions which bifurcates to the right of the
trivial one.}\\

\smallskip
 This lemma is a classical result of global bifurcation theory (for details see for example [C-R]). Indeed,
let us consider a positive T-periodic solution: if \ $T\neq T_k$\ , then the linearized associate
equation is non-singular.\\ We may also
deduce from bifurcation theorem, applied to the simple eigenvalues problem, that there is an
unique curve of non trivial solutions near the point \ $(T_k,\alpha).$ \  In fact, this uniqueness is
global. The trivial curve is \ $u_T \equiv\alpha$ .\\ Moreover, \
${\displaystyle ({du\over{dT}})_{T=T_k}}$ \ is an eigenvalue of the linearized associate equation.
According to the global bifurcation theory, we assert that the non trivial curves turn off on the
right of the singular solution \ $(T_k,u_{T_k}(t))$ \ . Consequently, when \ $T$ \ varies, two non
trivial curves never cross.\\

 To prove Theorem 1, we also need the following.\\

{\bf Lemma 4 }\qquad
{\it The minimal period  \ $T(c)$ \ of a (periodic, positive) solution\ $u_T$\ of the equation (1)
is a monotone increasing function of its energy, when \ $c \in [0,c_0]$}\\

\smallskip

 Chow-Wang [C-Wa] also have calculated the derivative of the period function\ $T(c)$\ for equation (2) and have found the
expression of the derivative
\\
$$ T'(c) = {1\over c} \int_a^b {{\phi ^2(w) - 2G(w)\phi '(w)}\over{\phi ^2(w)\sqrt{c - G(w)}}}dw.$$
Lemma 4  is a consequence of the following (see corollary (2-5) in [C-Wa]) \\

\bigskip
  {\bf Lemma 5 }\qquad
 {\it Consider the smooth function \ $g $\ such that \ $g(1) = 0 \ ,\ g'(1) > 0.$ \ 
Suppose that
 $$H(x) = g^2(x) - 2G(x)g'(x) +{{g''(1)}\over{3{g'}^2(1)}} g^3(x) > 0,$$ 
for all \ $x \in [a,b],$\ where\  $a < 1 < b$\ and \ $x \neq 1$ .\ Then, \qquad $T'(c) \geq
 0$ \qquad  for all  \ $c \in [0,c_0$ . \\ 
Moreover, if in addition 
$$ g''(x) \geq 0 \quad and \quad \Delta(x) = (x-1)\bigg[g'(x)g''(1) - g'(1)g''(x)\bigg] \geq  0
$$
then} \ $H(x) > 0$.

\bigskip 

 Notice that Lemma 4 implies Lemma 3 (see corollary (3-1) in [C-Wa]).\\

In order to apply the preceding lemma it is more convenient to 
make a change of variables
\
$h(t) =
\alpha f(\beta t)\ ,$
\ where 
$$ \alpha = \bigg({R\over{4(n-1) C}}\bigg)^{4/n} \qquad {\mbox and} \qquad \beta = \sqrt{{nC\over 4}}.$$
Note that the constant \ $C$ \ must be positive to ensure that the equation (1) to have a periodic
non constant solution. This change gives the equation
\begin{equation}
  f'' - f^{1-4/n} + f = 0.
\end{equation}
We can verify that this equation satisfies the Lemmata 2 and 3 given above. They will be used to
complete the analysis of the equation (1).\\ Indeed, we get the functions \\ $g(f) = f -
f^{1-{4\over n}} \quad;\quad g'(f) = 1 - (1-{4\over n})f^{-{4\over n}} \\ and \quad g''(f) =
{4\over n}(1-{4\over n})f^{-1-{4\over n}}$. \\ Notice that the hypothesis \quad $g''(f) \geq 0$\quad
is only satisfied in dimension \ $n \geq 4$.\\ Then, we calculate  $$\Delta(f) =
(f-1)[g''(1)g'(f) - g'(1)g''(f)].$$ \quad We get 
$$\Delta(f) = (f-1) {4\over n}(1-{4\over n}) [1 - f^{4\over n} + {4\over n}(f^{-{4\over n}} -
f^{-1-{4\over n}})].$$
It follows 
$$\Delta(f) = {4\over n}(1-{4\over n}) f(1-f^{-1})[1 - f^{4\over n} + {4\over n}f^{-{4\over n}}(1 -
f^{-1})],$$ 
which is obviously positive.\\ 
Concerning the case \ $n \leq 4$.\ Remark that for \ $n = 4$\ Equation (1) becomes linear. For \ $n=3$,\  we only notice that we get the implications $$g''(x) < 0 \Rightarrow  -2 G(x) g''(x) \geq 0 \Rightarrow H(x) > 0 $$ 
Hence, we are able to complete the proof of Theorem1. First, we deduce the increase of the period function depending on the energy of
the equation (1). This fact allows us to determine the lower bound of the number of these Derdzinski
metrics. We also remark, that the bifurcation points of the solution family are \ $$(T_k,u_k) \quad
,\ where
\quad u_k = ({(n-1)C\over{nR}})^{-n/4}\quad and \quad T_k = {{2\pi k}\over{\sqrt C}}.$$ 

It appears from the above analysis that, a non constant, periodic solution of the equation (1)
exists only if the circle length satisfies the following condition
\begin{equation}
  T > {{2\pi}\over{\sqrt C}}.
\end{equation}
Notice that, we get an infinity of solutions . All are obtained by rotation-translation of the variable.\\ Moreover, in the case
where \ $T$ \ satisfies the double inequality 
\begin{equation}
2\pi {{(k-1)}\over{\sqrt C}} < T \leq 2\pi{k\over{\sqrt C}} \quad k \ {\it is \  an \  integer} \  > 1, 
\end{equation}
 then the
equation (2) may admit at least \ k \ rotationnally invariant distinct solutions. Hence, we have improved the lemma 1 of
Derdzinski [D] (see also Chapter 16.33 of A. Besse [B]) \\

\section{ Pseudo-cylindric and Derdzinski metrics}
Notice that the Riemannian product \ $(S^1(T) \times N,dt^2 + h^{4/n}(t)g_0)$\ is conformally flat only if the factor \ $N$\ has constant sectional curvature.
 \\ Let  the Riemannian cylindric product \ $(S^1
\times S^{n-1},dt^2+d{\xi}^2), $ \ where \ $S^1$\ is the circle of length\ $T$\ and\ $(S^{n-1},d{\xi}^2)$\ 
is the standard sphere. A such metric has a parallel Ricci tensor. Moreover, we know the number of Yamabe metrics  is finite in the conformal class of the cylindric metric  \ $[dt^2+d{\xi}^2],$ \ (see [Ch], [M-P]).\\ We call {\it pseudo-cylindric metric} any non trivial Yamabe metric  $g_c$\ in \ $[dt^2+d{\xi}^2].$ (  \ $g_c$ \ is a Yamabe metric on a n-dimensional Riemannian manifold \ $(M,g)$\  means there is a \ $C^\infty $\ positive solution \ $u_c$ \  of a differential equation such that the metric \ $g_c = u_c^{4\over{n-2}} g $\  has a constant scalar curvature). \\  
For \ $k = 2$,\ there is a conformal diffeomorphism between \ $S^n -\{p_1,p_2\}$\ and  \ $(S^1
\times S^{n-1},dt^2+d{\xi}^2), $ \ where \ $S^1$\ is the circle of length\ $T$.\ The non trivial Yamabe metrics on \ $(S^1
\times S^{n-1},dt^2+d{\xi}^2), $ \ are called pseudo-cylindric metrics. There are metrics of the form \ $ g = u^{4\over{n-1}} (dt^2+d{\xi}^2)$ \ where the \ $C^\infty $\ function \ $u$ \ is a non constant positive solution of the Yamabe equation, {Ch]. \\   

We first remark there is a conformal diffeomorphism between the manifolds \quad
$I\!R^n -
\{0\}$
\quad and \quad $I\!R\times S^{n-1}$ \quad given by sending the point \ $x$ \ to  \ $(log\vert x
\vert, {x\over{\vert x \vert}})$. \  By using the stereographic projection, we see easily that the
manifold
\
$I\!R
\times S^{n-1}$ \ (which is the universal covering space of \ $S^1\times S^{n-1}$) \  is 
conformally equivalent to \ $S^n -\{0,\infty\}$ \ .\\ Notice that the manifold \ $S^n - (p,-p)$ \
with the standard induced metric, can be considered as the warped product \ $$]0,\pi[ \times
S^{n-1}\ ,
\quad\mbox{with allowed metric} \quad dt^2 + \sin^{2}t d\xi^2$$ It has been shown by
Cafarelli-Gidas-Spruck [C-G-S] ( using an Alexandrov reflection argument), that any solution of (2-1)
is in fact a spherically symmetric radial function (depending on geodesic distance from either \
$p$ \ or \
$-p$ \ ). Any solution of Equation (6) which gives a complete metric on the cylinder \ $\Re \times
S^{n-1}$ \ is of the form \ $u(t,\xi) = u(t)$ ,\ where \ $t \in I\!R $ \ and \ $\xi \in S^{n-1}$.\\
The background metric on the cylinder is the product \ $g_0 = dt^2 + d\xi^2 .$ \  For the
convenience, we assume the sphere radius equal to 1.\\ Therefore, the partial differential equation
 
\begin{equation}
 4{{n-1}\over{n-2}} \Delta_{g_0} u + R_{g_0} u - R_g
u^{{n+2}\over{n-2}} = 0,    
\end{equation} 

is
reduced to an ODE.\\ 
 The cylinder has scalar curvature \ $R(g_0) = (n-1)(n-2)$ \ and \ $R(u^{4\over{n-2}}g_0) =
n(n-1).$  Thus \ $u = u(t)$ \ satisfies  
\begin{equation}
 {d^2\over{dt^2}}u - {(n-2)^2\over 4}u + {n(n-2)\over 4}u^{{n+2}\over{n-2}} = 0.
\end{equation}
 It follows that a pseudo-cylindric metric (constant scalar curvature metric) on the
product \ $(S^1(T)\times S^{n-1} , g_0)$ \ corresponds to a T-periodic positive solution of (8) and
conversely. The analysis of this equation shows us, that it has only one center  \ $(\beta ,0)$ \ .
 This corresponding to the (trivial) constant solution  $$ \beta  = ({{n-2}\over{n}})^{{n-2}\over
4},$$
We proved the following, [Ch]

\bigskip

 {\bf Proposition } \quad {\it Consider the product manifold\quad $(S^1(T) \times S^{n-1},g_0)$\quad. Under
the condition  $$T(c) > T_1 = {{2\pi}\over{\sqrt{n-2}}},$$
  on the circle length, the Riemannian curvatures of
the associate pseudo-cylindric metrics\quad $g_c = {u_c}^{4\over{n-2}}g_0$ \quad are harmonic and their Ricci
tensors are non parallel.\\ Moreover, any pseudo-cylindric metric may be identified to a Derdzinski metric
up to a conformal transformation.}\\

\bigskip

Actually, any Derdzinski metric may be identified with a pseudo-cylindric metric up to conformal
diffeomorphism, let \ $F$ \ . Then the metrics are relied $$dt^2+{f^2}(t)d{\xi}^2 =
{F^*}({{u^j}_c}^{4\over{n-2}}(dt^2+d{\xi}^2)),$$ where \ ${{u^j}_c}$ \ are the pseudo-cylindric solutions
belonging to the (same) conformal class.\\
Indeed, for the metric \quad $dt^2 + f^2(t)d{\xi}^2$\quad  we can write $$ dt^2 + f^2(t)d{\xi}^2 = f^{2}(t)[({dt\over{f}})^2 +
d{\xi}^2].
$$After a change of variables and by using the conformal flatness of the product metric, we get $$dt^2 + f^2(t)g_0 = \phi^{2}(\theta)
[d\theta^2 + d{\xi}^2]$$ which is conformally flat. \\ This result is not surprising, because any manifold carying a warped metric
product
\quad
$(S^1 \times N,dt^2 + h^{4/n}(t)g_0)$\quad with harmonic curvature is not conformally flat
unless \ $(N,g_0)$\ is a space of constant curvature. This manifold must be locally conformally
equivalent to the trivial product \quad $S^1\times N.$ \\  This product will be conformally flat
only if
\ $N$\ has constant sectional curvature.\\

On the other hand, we remark that any warped metric  \ $dt^2 +
h^{4/n}(t)g_0$ \ defined by the lemma (4-4) on the product manifold \ $(S^1\times
S^{n-1},dt^2+d\xi^2) $\ is conformal to a Riemannian metric product $d\theta^2 + d\xi^2.$ \ Here \
$\theta $ \  is a new $S^1$-parametrisation with length \quad ${\displaystyle
\int_{S^1}{dt\over{h^{2/n}(t)}}}$ .\\
 Furthermore, we have seen in a previous paper (Chouikha [Ch]), there exists a analytic deformation
in the conformal class \ $[g_T]$\ , of any warped metric \ $g_T = dt^2 + f^2(t)d\xi^2$\ on the
manifold
\
$(S^1\times S^{n-1},dt^2+d\xi^2),\quad n=4 \ or \ 6$\ , and satisfying the length condition  $$
{\displaystyle T = \int_{S^1}{dt\over{f(t)}}}\ .$$ Notice that the product metric (under the
length condition) belongs to the conformal class\ $[g_T].$\  This analytic family of metrics depends
on two parameters \
$(g_{\alpha,\beta})$ \ . They all have a constant positive scalar curvature, and satisfy the
condition: \ $g_{\alpha,0}$ \ is the warped metric \ $g_T$\ .\\ Moreover, by using the same
argument of [Ch], we are able to extend the latter result to dimension 3. More precisely, a
metric on the euclidean space  \ $I\!R^3$ \ for which the rotation group \ $SO(3)$ \ acts by
isometries, is in fact a warped product as \ $g_T = dt^2 + f^2(t)d\xi^2$\ , where \ $d\xi^2$ \ is
the standard metric on the sphere \ $S^2$ \ . We can verify that its scalar curvature is $$R = {{2 -
2 {f'}^2 - 4 f f''}\over{f^2}}\ .$$ 

Moreover, we know that every conformally flat manifold\ $(M,g)$ \ admits a Codazzi tensor which is
not a constant multiple of the metric. Let \ $b$ \ a such symmetric 2-tensor field; suppose \ $b$
\ has exactly 2 distinct eigenvalues\ $\lambda \ , \ \mu$ , has constant trace\ $tr_g(b) = c$ \ and
is non parallel. Following [De], these conditions give locally  
\begin{equation}
M = I\times N \quad
\mbox{with  allowed  metric}
\quad g = dt^2 + e^{2\psi} g_N
\end{equation}
 and 
\begin{equation}
 b = \lambda dt^2 + \mu e^{2\psi} g_N, \quad
\lambda ={c\over n} + (1-n) c e^{-n\psi}, \quad \mu = {c\over n} + c e^{-n\psi}\ .
\end{equation}

Conversely, for
any such data and for an arbitrary function \ $\psi$ \quad on \quad $I$ ,\ (15) defines a
Riemannian manifold \ $(M,g)$ \ with Codazzi tensor \ $b$ \ of this type. If we assume that the
positive function \ $h = e^{{n\over 2}\psi}$ \ satisfies   Equation (1)
$$ h'' - {nR\over{4(n-1)}}h^{1-4/n} = -{n\over 4}C h \qquad \mbox{for some 
constant}
\qquad C > 0,$$ then the
Ricci tensor is precisely a Codazzi tensor and, thus it is non parallel.\\ Therefore, \ $(M,g)$ \ is
isometrically covered by \quad  $(\!R \times N,dt^2 + h^{4/n}(t)g_N)$  \quad where \ $(N,g_N)$ \
is an Einstein manifold with positive scalar curvature ( see [B] , § 16.33).

\vspace{1cm}

{\it {\bf REFERENCES}}\\

\bigskip

[Be]  A.L.Besse \ {\it Einstein Manifolds}\quad Ergebnisse der
Mathematik und ihrer Grenzgebiete, 3. Folge, Band 10. Springer. New-York, Berlin (1987).\\

[Ch]  R.Chouikha \ {\it M{\'e}triques conform{\`e}ment plates et {\'e}quations de Yang-Mills}\quad
C.R. Math. Reports. Acad Sc. of Canada, vol XIII (1991) 7-12.\\

[Ch1]  R.Chouikha \ {\it Famille de m{\'e}triques conform{\`e}ment plates et
r{\'e}gularit{\'e}}\quad C.R. Math. Reports. Acad Sc. of Canada, vol XIV (1992) 257-262.\\

[Ch2]  R. Chouikha \ {\it On the properties of certain metrics with constant scalar curvature}\quad
 Diff. Geom. and Appl., Sat. Conf. of ICM in Berlin, (1999), Brno, 39-46.\\ 

[Ch-W] R.Chouikha and F.Weissler \ {\it Monotonicity properties of the period function and the number of constant positive scalar curvature metrics on $S^1\times S^{n-1}$}\quad
Prepubl.  Math., Universite Paris-Nord, n° 94-14, (1994).\\

[C-Wa]  S.N.Chow and D.Wang \ {\it On the monotonicity of the period function of some second order
equations}\quad Casopi pro pestovani matematiky, 111 (1986), 14-25.\\

[C-R]  M.Crandall and P.Rabinowitz \ {\it Bifurcation from simple eigenvalues }\quad J. Functional
Analysis 8, (1971), 321-340.\\

[De]  A.Derdzinski \ {\it Classification of certain compact Riemannian manifolds with harmonic
curvature and non parallel Ricci tensor} \quad Math. Z., 172, (1980), 277-280.\\

[De1]  A. Derdzinski \ {\it On compact Riemannian manifolds with harmonic curvature } \quad Math.
Ann., vol 259, (1982), 145-152.\\

[S]  J.A. Shouten \ {\it Ricci Calculus} \quad Springer-Verlag, Berlin-Heidelberg, (1954).

\end{document}